# THE ASYMPTOTIC DISTRIBUTION OF A CLUSTER-INDEX FOR I.I.D. NORMAL RANDOM VARIABLES

By Yannis G. Yatracos

*National University of Singapore*

In a sample variance decomposition, with components functions of the sample's spacings, the largest component $\tilde{I}_n$ is used in cluster detection. It is shown for normal samples that the asymptotic distribution of $\tilde{I}_n$ is the Gumbel distribution.

**1. Introduction.** Clusters are nowadays data structures of considerable interest: Microarray data is used to attribute genes in clusters; gene expression is used to cluster tumors and identify similar types of cancer. Extreme value theory, in particular of sample spacings, has been used extensively in modeling phenomena. The extreme value $\tilde{I}_n$ of functions of spacings is introduced in Yatracos (2007) to detect data clusters from their one dimensional data projections; $n$ is the size of the data. In this work, the asymptotic distribution of $\tilde{I}_n$ is obtained for data from the normal distribution, and can be used to determine statistical significance of potential clusters.

Consider a sequence $X_1, \dots, X_n$ of independent identically distributed random variables with cumulative distribution function $F$. Let $X_{(i)}$ be the $i$th order statistic, $i = 1, \dots, n$. Define the spacing

$$S_i = X_{(i+1)} - X_{(i)}, \qquad i = 1, \dots, n-1,$$

the maximum spacing

$$M_n = \max\{S_i, i = 1, \dots, n-1\} = M_n^{(1)}$$

and the $k$th largest spacing $M_n^{(k)}, k = 1, \dots, n-1$,

The large sample behavior of $M_n$ and $M_n^{(k)}$, that is, their asymptotic distribution, large deviation properties and a.s. behavior has been studied for various choices of $F$ by Pyke (1965), Slud (1977/78), Devroye (1981, 1982, 1984), Deheuvels (1982, 1983, 1984, 1985) and other authors.

When $F = \Phi$, the cumulative distribution function of a standard normal $N(0,1)$ random variable, it is shown herein that the asymptotic distribution of

$$\tilde{I}_n = \max\{S_i\tilde{T}_i, i = 1, \ldots, n-1\}$$

is a standard Gumbel distribution;

$$\tilde{T}_i = \frac{i(n-i)}{n^2} \frac{(\bar{X}_{[i+1,n]} - \bar{X}_{[1,i]})}{\sum_{i=1}^n (X_i - \bar{X})^2/n}, \qquad i = 1, \ldots, n-1,$$

$\bar{X}_{[i,j]}$ is the average of the order statistics from $i$ to $j, i < j$.

Results in Deheuvels (1985) are crucial to obtain the result.

$S_i\tilde{T}_i$ is the $i$th component in a standardized sample variance decomposition, $i = 1, \ldots, n-1$ [Yatracos (1998)]. The largest component in the decomposition, $\tilde{I}_n$, determines two least homogeneous sample clusters. For multivariate data, $\tilde{I}_n$ is used to determine two clusters with the least homogeneous one-dimensional data projection [Yatracos (2007)]. Significance with respect to the normal model is justified since for many high dimensional data sets to find unusual projections one should search for nonnormality [Diaconis and Freedman (1984)].

## 2. The sample variance decomposition and $\tilde{I}_n$.

Univariate observations $X_1, \ldots, X_n$ are usually separated in two clusters by comparing the standardized difference of the group averages $\bar{X}_{[1,i]}$ and $\bar{X}_{[i+1,n]}$, respectively, of the $i$ smaller observations $X_{(1)}, \ldots, X_{(i)}$ and of the $n-i$ larger observations $X_{(i+1)}, \ldots, X_{(n)}$, for $i = 1, \ldots, n-1$. The spacing $S_i$ between the two groups may vary and it is natural to be used in a dissimilarity measure. The product $\frac{i(n-i)}{n^2}(\bar{X}_{[i+1,n]} - \bar{X}_{[1,i]})S_i$ is related to the sample variance [Yatracos (1998)]

$$\frac{1}{n}\sum_{i=1}^n (X_i - \bar{X})^2 = \sum_{i=1}^{n-1} \frac{i(n-i)}{n^2}(\bar{X}_{[i+1,n]} - \bar{X}_{[1,i]})S_i$$

and measures between-groups variation. The standardized variance components

$$S_i\tilde{T}_i = \frac{i(n-i)}{n^2} \frac{(\bar{X}_{[i+1,n]} - \bar{X}_{[1,i]})S_i}{\sum_{i=1}^n (X_i - \bar{X})^2/n}, \qquad i = 1, \ldots, n-1$$

indicate the relative contribution of the groups $X_{(1)}, \ldots, X_{(i)}$ and $X_{(i+1)}, \ldots, X_{(n)}$ in the sample variability.

The statistic

$$\tilde{I}_n = \max\{S_i\tilde{T}_i, i = 1, \ldots, n-1\}$$

determines two potential clusters. When $\tilde{I}_n = S_j\tilde{T}_j$, these clusters are $\tilde{\mathcal{C}}_1 = \{X_{(1)}, \ldots, X_{(j)}\}$, $\tilde{\mathcal{C}}_2 = \{X_{(j+1)}, \ldots, X_{(n)}\}$ and the cluster separators are $\tilde{s}_1 = X_{(j)}, \tilde{s}_2 = X_{(j+1)}$.

## 3. The asymptotic distribution of $\tilde{I}_n$.

THEOREM 3.1. *Let $Z_1, \dots, Z_n$ be i.i.d. standard normal random variables, $x \in R$. Then it holds that*

$$\lim_{n\to+\infty} P[n\tilde{I}_n < x + \log n] = \exp\{-\exp\{-x\}\}. \tag{1}$$

The proof of Theorem 3.1 is based on the four lemmas that follow. It is enough to obtain the asymptotic distribution of

$$\max\left\{\frac{i(n-i)}{n^2}(\bar{Z}_{[i+1,n]} - \bar{Z}_{[i,n]})(Z_{(i+1)} - Z_{(i)})\right\}, \qquad i = 1, \dots, n-1.$$

Let $Z_{(i)}$ be the $i$th order statistic with density $g_i, i = 1, \dots, n$. Let $n_+$ (resp. $n_-$) be the number of positive (resp. negative) observations. Then it holds that $n_+ \sim n_- \sim \frac{n}{2}$ a.s.; $a_n \sim b_n$ denotes $\lim_{n\to\infty} \frac{a_n}{b_n} = 1$. Due to the symmetry of the normal distribution, without loss of generality, the lemmas are proved for the *positive* observations $0 < Z_{(n/2+l_n)}, \dots, Z_{(n)}, Z_{(n/2+l_n-1)} < 0, l_n = o(n)$ can take either positive or negative values. One may think of the arguments as conditional on the value of $n_+$.

LEMMA 3.1. *For $\varepsilon > 0$ and $i = \frac{n}{2} + l_n, \dots, n$, it holds that*

$$P[Z_{(i)}(Z_{(i+1)} - Z_{(i)}) > \varepsilon] \le (1 - \varepsilon e^{-1.5\varepsilon})^{n-i}. \tag{2}$$

PROOF. Recall that for any $x > 0$ it holds

$$\Phi\left(x + \frac{\varepsilon}{x}\right) - \Phi(x) \ge \frac{\varepsilon\phi(x + \varepsilon/x)}{x}, \tag{3}$$

$$\frac{x}{1+x^2}\phi(x) < 1 - \Phi(x) < \frac{\phi(x)}{x} \tag{4}$$

[Chow and Teicher (1988), page 49],

and thus

$$\frac{\Phi(x+\varepsilon/x) - \Phi(x)}{1 - \Phi(x)} \ge \frac{\varepsilon\phi(x+\varepsilon/x)}{\phi(x)} = \varepsilon e^{-0.5\varepsilon^2/x^2 - \varepsilon}. \tag{5}$$

The Markovian property of $Z_{(1)}, \dots, Z_{(n)}$ implies that given $Z_{(i)} = z$, the r.v.'s $Z_{(i+1)}, \dots, Z_{(n)}$ form a sample from a standard normal distribution truncated at $z$ and, therefore,

$$\begin{aligned} &EP(Z_{(i)}(Z_{(i+1)} - Z_{(i)}) > \varepsilon | Z_{(i)} = x) \\ &\quad = \int_0^{+\infty} \left[\frac{1 - \Phi(x+\varepsilon/x)}{1 - \Phi(x)}\right]^{n-i} g_i(x)\,dx. \end{aligned} \tag{6}$$

For $0 < x \leq \sqrt{\varepsilon}$,

$$\left[\frac{1-\Phi(x+\varepsilon/x)}{1-\Phi(x)}\right]' = \frac{-\phi(x+\varepsilon/x)(1-\varepsilon/x^2)(1-\Phi(x))+\phi(x)(1-\Phi(x+\varepsilon/x))}{(1-\Phi(x))^2} > 0,$$

thus,

$$\frac{1-\Phi(x+\varepsilon/x)}{1-\Phi(x)} \leq \frac{1-\Phi(2\sqrt{\varepsilon})}{1-\Phi(\sqrt{\varepsilon})} = 1 - \frac{\Phi(2\sqrt{\varepsilon})-\Phi(\sqrt{\varepsilon})}{1-\Phi(\sqrt{\varepsilon})} \leq 1-\varepsilon e^{-1.5\varepsilon}. \tag{7}$$

The last inequality follows from (5) with $x = \sqrt{\varepsilon}$.

For $x > \sqrt{\varepsilon}$, it holds that

$$e^{-0.5\varepsilon} < e^{-0.5\varepsilon^2/x^2}. \tag{8}$$

From (5) and (8), it follows that

$$\begin{aligned}\frac{1-\Phi(x+\varepsilon/x)}{1-\Phi(x)} &= 1 - \frac{\Phi(x+\varepsilon/x)-\Phi(x)}{1-\Phi(x)} \\ &\leq 1-\varepsilon e^{-0.5\varepsilon^2/x^2-\varepsilon} < 1-\varepsilon e^{-1.5\varepsilon}.\end{aligned} \tag{9}$$

Inequality (2) follows from (6), (7) and (9). □

Lemma 3.2. *Let* $R_i = \frac{\phi(Z_{(i)})}{\phi(Z_{(i)}+T_i(Z_{(i+1)}-Z_{(i)}))}$, $T_i$ *is a random variable whose existence is guaranteed by Taylor's theorem,* $0 < T_i < 1, i = \frac{n}{2}+l_n,\dots,n-1, k_n \sim (\log n)^{1+\zeta}, 0<\zeta<2, m_n \to +\infty$ *as* $n \to +\infty$. *Then for small* $\varepsilon > 0$ *as* $n \to +\infty$, *it holds that*

$$P\left[\sup\left\{R_i, i = \frac{n}{2}+l_n,\dots,n-k_n\right\} > 1+\varepsilon\right] \to 0, \tag{10}$$

$$P[\sup\{R_i, i = n-k_n+1,\dots,n-1\} > m_n] \to 0. \tag{11}$$

Proof. For $i = \frac{n}{2}+l_n,\dots,n-1$, it holds that

$$1 \leq R_i \leq e^{0.5(Z_{(i+1)}-Z_{(i)})^2+Z_{(i)}(Z_{(i+1)}-Z_{(i)})}. \tag{12}$$

In Deheuvels (1985), it is shown that for $\eta > 0$,

$$P[\sqrt{2\log n}\max\{Z_{(i+1)}-Z_{(i)}; i=1,\dots,n-1\} > (1+\eta\log\log n) \text{ i.o.}] = 0;$$

i.o. denotes "infinitely often." Thus, from (12), to prove (10) and (11), it is enough to prove respectively that as $n \to +\infty$,

$$P\left[\sup\left\{Z_{(i)}(Z_{(i+1)}-Z_{(i)}), i = \frac{n}{2}+l_n,\dots,n-k_n\right\} > \varepsilon\right] \to 0,$$

$$P[\sup\{Z_{(i)}(Z_{(i+1)}-Z_{(i)}), i = n-k_n+1,\dots,n-1\} > \log m_n] \to 0.$$

From (2), it follows that

$$P\left[\sup\left\{Z_{(i)}(Z_{(i+1)}-Z_{(i)}), i=\frac{n}{2}+l_n,\dots,n-k_n\right\}>\varepsilon\right]$$
$$\leq\left(\frac{n}{2}-k_n-l_n+1\right)(1-\varepsilon e^{-1.5\varepsilon})^{k_n}\to 0$$

as $n\to+\infty$ since $l_n=o(n)$.[1]

For $\theta>0$, $Z_{(n)}\leq(1+\theta)\sqrt{2\log n}$ a.s. for $n\geq n(\theta)$ and

$$P[\sup\{Z_{(i)}(Z_{(i+1)}-Z_{(i)}), i=n-k_n+1,\dots,n-1\}>\log m_n]$$
$$\leq P[(1+\theta)\sqrt{2\log n}$$
$$\times\sup\{(Z_{(i+1)}-Z_{(i)}), i=n-k_n+1,\dots,n-1\}>\log m_n].$$

From Lemma 6 in Deheuvels (1985), the $K_n=[(\log n)^3]$ largest order statistics generate spacings which are uniformly close to $(2\log n)^{-1/2}E_j/j, j=1,\dots,K_n$, where $\{E_j, j=1,\dots,K_n\}$ are i.i.d. exponential r.v.'s with mean 1. Thus, it holds that

$$P\left[\sqrt{2\log n}(Z_{(j+1)}-Z_{(j)})>\frac{\log m_n}{1+\theta}\right]$$
$$\sim P\left[E_j>\frac{j\log m_n}{1+\theta}\right]=e^{-j\log m_n/(1+\theta)},\qquad j=1,\dots,K_n$$

and

$$P\left[\sqrt{2\log n}\sup\{(Z_{(i+1)}-Z_{(i)}), i=n-k_n+1,\dots,n-1\}>\frac{\log m_n}{1+\theta}\right]$$
$$\leq\sum_{j=1}^{k_n-1}e^{-j\log m_n/(1+\theta)}=e^{-\log m_n/(1+\theta)}\frac{1-e^{-(k_n-1)\log m_n/(1+\theta)}}{1-e^{-\log m_n/(1+\theta)}}$$
$$\sim e^{-\log m_n/(1+\theta)}\longrightarrow 0$$

as $n\to+\infty$. □

Lemma 3.3. *For any real $x$, as $n\to+\infty$,*

$$P\left[\bigcap_{i=n/2+l_n}^{n-1}\{(n-i)(Z_{(i+1)}-Z_{(i)})E(\bar{Z}_{[i+1,n]}|Z_{(i)})\leq x+\log n\}\right]\to e^{-0.5e^{-x}}.$$

Proof. Let $k_n\sim(\log n)^{1+\zeta}, 0<\zeta<2$, $m_n\sim\log\log n$,

$$A_i=\{(n-i)(Z_{(i+1)}-Z_{(i)})E(\bar{Z}_{[i+1,n]}|Z_{(i)})\leq x+\log n\}\tag{13}$$

[1]The result follows also from Deheuvels (1985) with $k_n=[(\log n)^3]$.

and $A_i^c$ its complement, $i = \frac{n}{2} + l_n, \dots, n-1$.

$$P\left[\bigcap_{i=n/2+l_n}^{n-k_n} A_i \cap \left(\bigcap_{i=n-k_n+1}^{n-1} A_i\right)\right]$$
$$= P\left[\bigcap_{i=n/2+l_n}^{n-k_n} A_i\right] - P\left[\bigcap_{i=n/2+l_n}^{n-k_n} A_i \cap \left(\bigcup_{i=n-k_n+1}^{n-1} A_i^c\right)\right]$$

and it is enough to show that as $n \to +\infty$

$$P\left[\bigcap_{i=n/2+l_n}^{n-k_n} A_i\right] \to e^{-0.5e^{-x}}$$

and

$$P\left[\bigcup_{i=n-k_n+1}^{n-1} A_i^c\right] \to 0.$$

Let $U_{(i)}$ be the $i$th order statistic from $n$ i.i.d. uniform r.v.'s on $(0,1)$. Then it holds that

$$Z_{(i+1)} - Z_{(i)} = \Phi^{-1}(U_{(i+1)}) - \Phi^{-1}(U_{(i)}) = \frac{(U_{(i+1)} - U_{(i)})}{\phi(\tilde{Z}_i)}, \tag{14}$$

with $\tilde{Z}_i = Z_{(i)} + T_i(Z_{(i+1)} - Z_{(i)}), T_i$ is a random variable whose existence is guaranteed by Taylor's theorem, $0 < T_i < 1$.

Given $Z_{(i)} = z$, the r.v.'s $Z_{(i+1)}, \dots, Z_{(n)}$ form a sample from a standard normal distribution truncated at $z$ with mean $\frac{\phi(z)}{1-\Phi(z)}$ and variance $1 + \frac{z\phi(z)}{1-\Phi(z)} - [\frac{\phi(z)}{1-\Phi(z)}]^2$. Thus, it holds that

$$E(\bar{Z}_{[i+1,n]}|Z_{(i)}) = \frac{\phi(Z_{(i)})}{1 - \Phi(Z_{(i)})}$$

and that

$$(n-i)(Z_{(i+1)} - Z_{(i)})E(\bar{Z}_{[i+1,n]}|Z_{(i)}) = (n-i)\frac{(U_{(i+1)} - U_{(i)})}{\phi(\tilde{Z}_i)}\frac{\phi(Z_{(i)})}{1 - \Phi(Z_{(i)})}$$
$$= (n-i)\left(1 - \frac{V_{(n-i)}}{V_{(n-i+1)}}\right)R_i,$$

$R_i$ is defined as in Lemma 3.2, $V_{(n-i)} = 1 - \Phi(Z_{(i)})$ is the $(n-i)$th order statistic from i.i.d. uniform r.v.'s on $(0,1)$, and

$$\left(\frac{V_{(n-i)}}{V_{(n-i+1)}}\right)^{n-i}, \qquad i = 1, \dots, n-1$$

are i.i.d. uniform random variables on $(0,1)$ (see, e.g., David and Nagaraja (2003)).

From Lemma 3.2, using $D_n = x + \log n$, it holds that

$$P\left[\bigcap_{i=n/2+l_n}^{n-k_n} A_i\right] \sim P\left[\bigcap_{i=n/2+l_n}^{n-k_n} \left\{(n-i)\left(1-\frac{V_{(n-i)}}{V_{(n-i+1)}}\right) \le D_n\right\}\right]$$
$$= \prod_{i=n/2+l_n}^{n-k_n} \left[1-\left(1-\frac{D_n}{n-i}\right)^{n-i}\right]$$
$$\sim \prod_{i=n/2+l_n}^{n-k_n} (1-e^{-D_n}) = (1-e^{-\log n - x})^{n/2-k_n-l_n+1} \sim e^{-0.5e^{-x}},$$

since $l_n = o(n)$.

From (11), it also holds that

$$P\left[\bigcup_{i=n-k_n+1}^{n-1} A_i^c\right] \le \sum_{i=n-k_n+1}^{n-1} P\left[(n-i)m_n\left(1-\frac{V_{(n-i)}}{V_{(n-i+1)}}\right) > D_n\right]$$
$$= \sum_{i=n-k_n+1}^{n-1} \left(1-\frac{D_n}{m_n(n-i)}\right)^{n-i}$$
$$\le k_n e^{-(x+\log n)/m_n} \to 0$$

as $n \to +\infty$. $\square$

Lemma 3.4. *Let* $k_n \sim (\log n)^{1+\zeta}, 0<\zeta<2$. *As* $n\to+\infty$:

(a) $\sup\{(n-i)(Z_{(i+1)}-Z_{(i)})|\bar{Z}|, i=\frac{n}{2}+l_n,\ldots,n-1\} \to 0$,

(b) $\sup\{(n-i)(Z_{(i+1)}-Z_{(i)})|[\bar{Z}_{[i+1,n]} - E(\bar{Z}_{[i+1,n]}|Z_{(i)})]|, i=\frac{n}{2}+l_n,\ldots, n-k_n\} \to 0$,

(c) $\sup\{(n-i)(Z_{(i+1)}-Z_{(i)})|[\bar{Z}_{[i+1,n]} - E(\bar{Z}_{[i+1,n]}|Z_{(i)})]|, i=n-k_n+1,\ldots,n-1\} \to 0$,

*all in probability.*

Proof. (a) (4) implies that $\frac{1-\Phi(x)}{\phi(x)}$ is decreasing for $x\ge 0$. Thus, it holds that

$$\sup\left\{(n-i)(Z_{(i+1)}-Z_{(i)}), i=\frac{n}{2}+l_n,\ldots,n-1\right\}|\bar{Z}|$$
$$\le \sup\left\{(n-i)(Z_{(i+1)}-Z_{(i)})E(\bar{Z}_{[i+1,n]}|Z_{(i)}), i=\frac{n}{2}+l_n,\ldots,n-1\right\}$$
$$\times \frac{1-\Phi(0)}{\phi(0)}|\bar{Z}| \to 0$$

in probability as $n \to +\infty$; use Lemma 3.3 and limit theorems for $\bar{Z}$.

(b) Conditionally on $Z_{(i)}$, let $\sigma^2_{Z_{(i)}} = 1 + \frac{Z_{(i)}\phi(Z_{(i)})}{1-\Phi(Z_{(i)})} - [\frac{\phi(Z_{(i)})}{1-\Phi(Z_{(i)})}]^2 < 1$, denote the variance of $Z_{(i+1)}, \dots, Z_{(n)}, i = 0.5n + l_n, \dots, n - k_n$. Then it holds that

$$\sup\left\{(n-i)(Z_{(i+1)} - Z_{(i)})|\bar{Z}_{[i+1,n]} - E(\bar{Z}_{[i+1,n]}|Z_{(i)})|, i = \frac{n}{2} + l_n, \dots, n - k_n\right\}$$
$$\leq \frac{1-\Phi(0)}{\phi(0)} \sup\left\{(n-i)(Z_{(i+1)} - Z_{(i)})E(\bar{Z}_{[i+1,n]}|Z_{(i)}) \times \frac{|\bar{Z}_{[i+1,n]} - E(\bar{Z}_{[i+1,n]}|Z_{(i)})|}{\sigma_{Z_{(i)}}}, i = \frac{n}{2} + l_n, \dots, n - k_n\right\}.$$

Let $S_i = \sum_{j=i+1}^n Z_{(j)} - E(\sum_{j=i+1}^n Z_{(j)}|Z_{(i)}), i = 0.5n + l_n, \dots, n - k_n$. For $\delta > 0$, it holds that

$$P\left[(n-i)^{0.1}\frac{|\bar{Z}_{[i+1,n]} - E(\bar{Z}_{[i+1,n]}|Z_{(i)})|}{\sigma_{Z_{(i)}}} > \delta\right]$$
$$= EP\left[\frac{|S_i|}{\sigma_{Z_{(i)}}\sqrt{n-i}} > \delta(n-i)^{0.4}|Z_{(i)} = z\right]$$
$$\leq \left|EP\left[\frac{|S_i|}{\sigma_z\sqrt{n-i}} > \delta(n-i)^{0.4}|Z_{(i)} = z\right] - P[|Z| > \delta(n-i)^{0.4}]\right| + P[|Z| > \delta(n-i)^{0.4}]$$
$$\leq \frac{2C_1C_U}{\sqrt{n-i}[1 + (n-i)^{0.8}\delta^2]} + C_2\frac{e^{-0.5\delta^2(n-i)^{0.8}}}{\delta(n-i)^{0.4}},$$

$C_U$ is the universal constant in the Berry–Esseen bound [see, e.g., Serfling (1980) or Ibragimov and Linnik (1971)], $C_2$ is positive constant. The Markovian property of $Z_{(1)}, \dots, Z_{(n)}$ implies that given $Z_{(i)} = z$, the r.v.'s $Z_{(i+1)}, \dots, Z_{(n)}$ form a sample from a standard normal distribution truncated at $z$, therefore,

$$\sup\left\{\frac{E(|Z_{(j+1)} - E(Z_{(j+1)}|Z_{(i)} = z)|^3|Z_{(i)} = z)}{\sigma_z^{3/2}}, z > 0, j = i, \dots, n-1\right\}$$

in the bound is replaced by its equivalent for the sample

$$\sup\left\{\frac{E[|Z_1 - E(Z_1|Z_1 > z)|^3|Z_1 > z]}{\sigma_z^{3/2}}, z > 0\right\} = C_1.$$

$C_1$ is bounded since:

(i) for $z>0$ large, (4) implies $z \sim \frac{\phi(z)}{1-\Phi(z)} = E(Z_1|Z_1>z)$ and

$$\frac{E[|Z_1 - E(Z_1|Z_1>z)|^3|Z_1>z]}{\sigma_z^{3/2}} \approx \frac{E[((Z_1 - E(Z_1|Z_1>z))^3)^+|Z_1>z]}{\sigma_z^{3/2}} \approx \frac{E[(Z_1 - E(Z_1|Z_1>z))^3|Z_1>z]}{\sigma_z^{3/2}},$$

where $a^+ = \max(a,0)$,

(ii) $\lim_{z\to+\infty} \frac{E[(Z_1-E(Z_1|Z_1>z))^3|Z_1>z]}{\sigma_z^{3/2}} = 2$ [Nariaki and Akihide (1985)],

(iii) $\frac{E[(Z_1-E(Z_1|Z_1>z))^3|Z_1>z]}{\sigma_z^{3/2}}$ is continuous function in $z$ and, therefore, achieves its maximum in any compact $[0,M], M>0$, and in particular for $M$ large.

Thus, as $n \longrightarrow +\infty$, it holds that

$$\sum_{i=0.5n+l_n}^{n-k_n} P\left[(n-i)^{0.1}\frac{|\bar{Z}_{[i+1,n]} - E(\bar{Z}_{[i+1,n]}|Z_{(i)})|}{\sigma_{Z_{(i)}}} > \delta\right] \leq \sum_{i=0.5n+l_n}^{n-k_n} \frac{2C_1C_U}{(n-i)^{1.3}} + C_2 \sum_{i=0.5n+l_n}^{n-k_n} \frac{e^{-0.5\delta^2(n-i)^{0.8}}}{(n-i)^{0.4}} \sim C_3\left(\frac{1}{(\log n)^{0.3(1+\zeta)}} - \frac{1}{(0.5n-l_n)^{0.3}}\right) \longrightarrow 0;$$

$C_3$ is a positive positive constant, $l_n = o(n)$. Let

$$G_i = \left\{(n-i)^{0.1}\frac{|\bar{Z}_{[i+1,n]} - E(\bar{Z}_{[i+1,n]}|Z_{(i)})|}{\sigma_{Z_{(i)}}} \leq \delta\right\}, \qquad i = 0.5n+l_n, \ldots, n-k_n.$$

Using relations in the proof of Lemma 3.3, it follows that

$$\sum_{i=0.5n+l_n}^{n-k_n} P\left[\left\{(n-i)(Z_{(i+1)} - Z_{(i)}) \times E(\bar{Z}_{[i+1,n]}|Z_{(i)})\frac{|\bar{Z}_{[i+1,n]} - E(\bar{Z}_{[i+1,n]}|Z_{(i)})|}{\sigma_{Z_{(i)}}} > \varepsilon\right\} \cap G_i\right]$$

$$\leq \sum_{i=0.5n+l_n}^{n-k_n} P\left[(n-i)(Z_{(i+1)} - Z_{(i)})E(\bar{Z}_{[i+1,n]}|Z_{(i)}) > \frac{\varepsilon}{\delta}(n-i)^{0.1}\right]$$

$$\leq \sum_{i=0.5n+l_n}^{n-k_n} e^{-\varepsilon/\delta(n-i)^{0.1}} \sim \int_{0.5n+l_n}^{n-(\log n)^{1+\zeta}} e^{-\varepsilon/\delta(n-x)^{0.1}}\,dx$$

$$= \sum_{k=1}^{10} C_k^* y^k e^{-cy}]_{(\log n)^{0.1(1+\zeta)}}^{(0.5n-l_n)^{0.1}} + C_{11}^* \int_{(\log n)^{0.1(1+\zeta)}}^{(0.5n-l_n)^{0.1}} e^{-cy}\, dy \to 0$$

as $n \to +\infty$; $c = \frac{\varepsilon}{\delta}, l_n = o(n), C_k^*$ is a constant, $k = 1, \ldots, 11$.

(c) From Deheuvels (1985), for $i = n - k_n + 1, \ldots, n-1$, it holds that

$$\sqrt{2\log n}(n-i)(Z_{(i+1)} - Z_{(i)}) \sim E_{n-i};$$

$E_j, j = 1, \ldots, k_n - 1$, are i.i.d. exponential with mean 1. From (4), it also holds that

$$\begin{aligned}
&\bar{Z}_{[i+1,n]} - E(\bar{Z}_{[i+1,n]}|Z_{(i)}) \\
&\quad\sim \frac{(Z_{(i+1)} - Z_{(i)}) + (Z_{(i+2)} - Z_{(i)}) + \cdots + (Z_{(n)} - Z_{(i)})}{n-i} \\
&\quad= ((n-i)(Z_{(i+1)} - Z_{(i)}) \\
&\qquad + (n-i-1)(Z_{(i+2)} - Z_{(i+1)}) + \cdots + (Z_{(n)} - Z_{(n-1)}))(n-i)^{-1} \\
&\quad\sim \frac{E_1 + \cdots + E_{n-i}}{\sqrt{2\log n}(n-i)} \le \frac{\sup\{E_j, j = 1, \ldots, n-i\}}{\sqrt{2\log n}} \\
&\quad\le \frac{\sup\{E_j, j = 1, \ldots, k_n - 1\}}{\sqrt{2\log n}}.
\end{aligned}$$

Thus,

$$\begin{aligned}
&P[\sup\{(n-i)(Z_{(i+1)} - Z_{(i)})|\bar{Z}_{[i+1,n]} - E(\bar{Z}_{[i+1,n]}|Z_{(i)})|, \\
&\qquad\qquad i = n - k_n + 1, \ldots, n-1\} > \varepsilon] \\
&\quad\le P\left[\sup\left\{\frac{E_{n-i}}{\sqrt{2\log n}} \frac{\sup\{E_j, j = 1, \ldots, k_n - 1\}}{\sqrt{2\log n}},\right.\right. \\
&\qquad\qquad \left.\left. i = n - k_n + 1, \ldots, n-1\right\} > \varepsilon\right] \\
&\quad= P[\sup\{E_j, j = 1, \ldots, k_n - 1\} > \sqrt{\varepsilon}\sqrt{2\log n}] \\
&\quad= 1 - (1 - e^{-\sqrt{\varepsilon}\sqrt{2\log n}})^{k_n - 1} \to 0,
\end{aligned}$$

as $n \to +\infty$. □

Proof of Theorem 3.1. Conditionally on the value of $n_+$, from the definition of $l_n$ (before Lemma 3.1) it holds that $Z_{(n/2+l_n-1)} < 0 < Z_{(n/2+l_n)}$, and assume without loss of generality that $\phi(Z_{(n/2+l_n-1)}) > \phi(Z_{(n/2+l_n)})$. Note also that it holds

$$\begin{aligned}
\frac{i(n-i)}{n^2}(\bar{Z}_{[i+1,n]} - \bar{Z}_{[1,i]}) &= \frac{n-i}{n}(\bar{Z}_{[i+1,n]} - \bar{Z}) \\
&= \frac{i}{n}(\bar{Z} - \bar{Z}_{[1,i]}), \qquad i = 1, \ldots, n-1.
\end{aligned}$$

Let $A_i$ be as in (13), $B_i = \{i(Z_{(i+1)} - Z_{(i)})(-1)E(\bar{Z}_{[1,i]}|Z_{(i+1)}) \leq x + \log n\}$, $i = 1, \dots, \frac{n}{2}$. From Lemma 3.3 and its proof, it follows directly for the $A$'s, and by symmetry for the $B$'s that

$$P\left[\bigcap_{i=k_n+1}^{n/2} B_i \bigcap_{i=n/2+1}^{n-k_n} A_i\right] \sim P\left[\bigcap_{i=k_n+1}^{n/2} B_i\right] P\left[\bigcap_{i=n/2+1}^{n-k_n} A_i\right] \sim e^{-e^{-x}},$$

and as $n \to +\infty$,

$$P\left[\bigcup_{i=1}^{k_n} B_i^c \bigcup_{n-k_n+1}^{n-1} A_i^c\right] \to 0.$$

The proof is completed using Lemma 3.4. $\square$

**Acknowledgments.** Many thanks are due to an anonymous referee for a very careful and thorough review, and for useful suggestions that served to improve the presentation of the paper. Thanks are also due to Professor Prabir Burman for useful suggestions in an earlier version of this paper.

Department of Statistics and Applied Probability
The National University of Singapore
6 Science Drive 2
Singapore 117546
Republic of Singapore
E-mail: yatracos@stat.nus.edu.sg
stayy@nus.edu.sg